# Bootstrap of means under stratified sampling

**Odile Pons**

*INRA, Mathématiques,*
*78352 Jouy en Josas cedex, France*
*e-mail:* `Odile.Pons@jouy.inra.fr`

**Abstract:** In a two-stage cluster sampling procedure, $n$ random populations are drawn independently from independent populations and a subsample of observations is taken in each of them. The estimator of the general mean of the observed variables is asymptotically Gaussian and the asymptotic distributions of several bootstrap versions of the normalized and studentized statistics are studied. A weighted population resampling provides a good approximation and its accuracy depends on the convergence rate of the sample size of the populations.

## Contents

## 1. Introduction

In a classical stratified sampling, a population is split into a fixed number $L$ of strata and in each of them a subsample is observed. Various bootstrap estimators of the mean have been studied as the strata sizes are finite [16, 18] or infinite [2, 4]. In a two-stage cluster sampling the population consists of clusters of units, first some clusters are sampled and then units are sampled within the selected clusters. Several two-stage boostrap methods have also been studied in the case of a finite number of finite clusters drawn without replacement [16, 18]. Here, we consider a two-stage cluster sampling with an infinite number of clusters. The general setting is the following: A general population is subdivided into a large number of independent populations (or clusters) which cannot be all observed. We consider the problem of estimating some parameters of the distribution function of a variable $X$ on individuals in the general population and we assume that the realizations of this variable in the different populations are independent and identically distributed. The observations of $X$ on individuals within the same population are supposed to be independent and identically distributed conditionally on the selected population.

Here the parameters of interest are the means of $X$, the variance of the means of $X$ in the different sub-populations and the variance of $X$ and we focus on the behavior of bootstrap estimators. The normalized and the studentized estimators of the mean are asymptotic Gaussian as the number of sampled populations tends to infinity and the variances are consistently estimated. Three bootstrap procedures are studied,

B1 Sampling the populations with replacement in the set of the observed populations, then taking the original observed data of these sampled populations,

B2 Taking each of the observed populations and sampling the individuals with replacement in the set of the observed data within each population,

B3 A two-stage bootstrap cluster sampling that is a combination of the first two procedures. Here the populations are sampled with replacement from the set of the observed populations, then individuals are sampled with replacement within each population.

An application and simulation of the B1 sampling for a large number of populations was presented in [11, 12] for specific parameters in forestry. Bootstrap sampling with a large number of dependent variables looks like a stratified sampling, where individuals and variables may be drawn in the bootstrap procedure. Only the bootstrap resampling of individuals is relevant due to the loss of dependence between dependent variables of an individual when they are not jointly resampled, as in [13, 14, 15] with an application to medical data.

The asymptotic properties of the mean bootstrap estimator in the classical i.i.d. bootstrap were studied by many authors (e.g.[1, 2, 3, 7, 17]). They have been extended to some cases of independent but non identically distributed variables and to functional results in [2, 4, 6, 9, 18]. We prove properties for some of the three bootstrap procedures when the number of strata and the sub-sampling sizes tend to infinity, under a condition for their respective rate. The results are quite different from those of [1, 2, 16, 18] for means of fixed sub-populations and with resampling of all the populations (B2 sampling scheme), equivalent results are established in proposition 2.2 for the mean of the sub-population means. Though the bootstrap resampling scheme B1 is consistent for the mean parameters $\mu$ and $\mu_k$ of the general population and sub-populations, it is not for the variances of their estimators $\widehat{\mu}$ and $\widehat{\mu}_k$ and other bootstrap estimators of the variances are necessary to obtain asymptotically normal estimators. With random populations, the bootstrap sampling scheme B2 gives consistent and asymptotically normal estimators for all the parameters with weighted bootstrap. Under a weighted bootstrap resampling scheme B1, the estimator of the global mean and its variance achieve the usual properties. Under sampling scheme B3, the bootstrap estimators of $\mu$ and $\mu_k$ and the variances are not consistent. The rate of the bootstrap approximation is given for the weighted resampling scheme B1. As the variance of the estimators splits into within and between-population variances of different orders, the asymptotic results are quite different from the usual results.

## 2. Empirical estimation in a subdivided population

A variable $X$ is observed according to the following two-stage cluster sampling: a sample of $K$ populations is selected among a large number $L$ of populations and, in the $k$-th population of unknown large size, we consider $n_k$ independent observations of the variable of interest, $X_{ki}, i \leq n_k, k \leq K$. Let $N = \sum_{k \leq K} n_k$ be the total number of observations, we suppose that $N$ and $n_k$, for each $k$, increase with $K$. The estimators of the global population have then equivalently indexed by $N = N_K$ or $K$.

Conditionally on the $k$-th population sampling, the variables $X_{ki}, i \leq n_k$, have the distribution function $F_k$ and the distribution of the variables $X_{ki}$ in the general population is $F$ including the distribution functions of the $L$ subpopulations. Denote by $E$ the expectation with respect to the sampling of the populations and by $E_k$ the conditional expectation in the $k$-th population, $\mu_k$ the conditional mean of the $X_{ki}$'s in the $k$-th population $\mu = EX_{ki} = E\mu_k$.

If $E|X|^2 < \infty$, the variance of $X$ is denoted $V$ in the general population and the random variance of the $X_{ki}$'s conditionally on the sampling of the $k$-th population is $V_k = E_k(X_{ki} - \mu_k)^2$. We also denote by $\gamma = E(\mu_k - \mu)^2$ the variance of the independent random variables $\mu_k$, and $\sigma^2 = EV_k$. The variance of the $X_{ki}$'s is

$$V = \sigma^2 + \gamma. \tag{1}$$

The variables $X_{ki}$ and $X_{kj}$, $i \neq j$, are independent within the $k$-th population but they are dependent in the general population, with $Cov(X_{ki}, X_{kj}) = \gamma$. The variables $X_{ki}$ and $X_{lj}$, $k \neq l$, are independent for any $i$ and $j$. A random effect linear model having a nested error structure can be used to describe the data,

$$X_{ki} = \mu + a_k + u_{ki}, i = 1, \ldots, n_k, k = 1 \ldots, K,$$

with $a_k = \mu_k - \mu$ and $u_{ki} = X_{ki} - \mu_k$, where $E_k u_{ki} = 0$, $Var_k u_{ki} = V_k$ and $E_k(u_{ki}u_{kj}) = 0$ for $i \neq j$, $Ea_k = 0, Var a_k = \gamma$ and $E(a_k a_l) = 0$ for $k \neq l$. Within the $k$-th population, $\mu_k$ and $V_k$ are unbiasedly estimated by

$$\begin{aligned} \widehat{\mu}_k &= \frac{1}{n_k} \sum_{i \leq n_k} X_{ki}, \\ \widehat{V}_k &= \frac{1}{n_k - 1} \sum_{i \leq n_k} (X_{ki} - \widehat{\mu}_k)^2. \end{aligned} \tag{2}$$

For the global population, $\mu$ is estimated by the empirical mean $\widehat{\mu}_N$ of the $X$'s. Several boostrap estimators of $\mu$ will be proved to converge to the empirical estimator $\widehat{\mu}'_K$ of the mean of the intra-populations means. For the general mean, the estimators are denoted

$$\widehat{\mu}_N = \frac{1}{N} \sum_{k \leq K} \sum_{i \leq n_k} X_{ki} = \sum_{k \leq K} \frac{n_k}{N} \widehat{\mu}_k, \qquad (3)$$

$$\widehat{\mu}'_K = \frac{1}{K} \sum_{k \leq K} \widehat{\mu}_k. \qquad (4)$$

Note that with random populations, both $\widehat{\mu}_N$ and $\widehat{\mu}'_K$ have the expectation $\mu$ since $E\widehat{\mu}_k = \mu$ for every $k$. It would not be with fixed populations. We have $Var_k \widehat{\mu}_k = E_k(\widehat{\mu}_k - \mu_k)^2 = n_k^{-1} V_k$ and

$$Var \widehat{\mu}_k = E(\widehat{\mu}_k - \mu)^2 = \gamma + \frac{EV_k}{n_k}.$$

The variance terms $\sigma^2$, $V$ and $\gamma$ of the variable $X$ are estimated by

$$\widehat{\sigma}_N^2 = \frac{1}{N} \sum_{k \leq K} n_k \widehat{V}_k = \frac{1}{N} \sum_{k \leq K} \frac{n_k}{n_k - 1} \sum_{i \leq n_k} (X_{ki} - \widehat{\mu}_k)^2,$$

$$\widehat{\gamma}_N = \frac{K}{K-1} \sum_{k \leq K} (\widehat{\mu}_k - \widehat{\mu}'_K)^2 - \frac{1}{K} \sum_{k \leq K} \frac{1}{n_k - 1} \sum_{i \leq n_k} (X_{ki} - \widehat{\mu}_k)^2, \qquad (5)$$

$$\widehat{V}_N = \widehat{\gamma}_N + \widehat{\sigma}_N^2. \qquad (6)$$

The variances of the estimators are

$$S_N^2 := Var \widehat{\mu}_N = \frac{1}{N^2}(\sum_k n_k EV_k) + \gamma n^* = \frac{\sigma^2}{N} + \gamma n^*, \qquad (7)$$

$$S_K'^2 := Var \widehat{\mu}'_K = \frac{1}{K}(\gamma + \frac{1}{K} \sum_k \frac{EV_k}{n_k}) = \frac{1}{K}(\gamma + \frac{\sigma^2}{\tilde{n}}), \qquad (8)$$

where $EV_k$ satisfies (1), $\tilde{n}^{-1} = K^{-1} \sum_k n_k^{-1} < 1$ is the harmonic mean of the $n_k$'s and $n^* = N^{-2} \sum_k n_k^2 < N^{-1} \max_{k=1}^{K} n_k$ with $n^* = K^{-1}$ when $n_k$ constant for every $k$, then the estimators are identical. The variance of $\widehat{\mu}'_K$ is smaller than the variance of $\mu$ when $\gamma \sigma^{-2} > \{(\tilde{n}K)^{-1} - N^{-1}\}\{n^* - K^{-1}\}^{-1}$ which may happen even with $K = 2$. If $n_k = c_k K^\alpha$ with $0 < \alpha < 1$ and bounded constants $c_k$ for every $k$, then $S_N < S'_K$. The variance of both estimators split into a within-populations variance and a between-populations variance, $Var_{inter} \widehat{\mu}_N = N^{-1} \sigma^2$ and $Var_{inter} \widehat{\mu}_N = n^* \gamma$.

By (1)-(8) and with $n_k = O(K^\alpha)$, asymptotically unbiased estimators of $Var \widehat{\mu}'_K$ and $Var \widehat{\mu}_N$ are

$$\widehat{Var} \widehat{\mu}'_K = \frac{1}{K-1} \sum_{k \leq K} (\widehat{\mu}_k - \widehat{\mu}'_K)^2, \qquad (9)$$

$$\widehat{Var} \widehat{\mu}_N = n^* \widehat{\gamma}_N + \frac{\widehat{\sigma}_N^2}{N^2}. \qquad (10)$$

Estimators of the within and between-populations variances follow as

$$\widehat{V}ar_{inter}\widehat{\mu}_N = n^*\widehat{\gamma}_N, \quad \widehat{V}ar_{intra}\widehat{\mu}_N = \frac{1}{N^2}\sum_{k\leq K} n_k \widehat{V}_k,$$

$$\widehat{V}ar_{inter}\widehat{\mu}'_K = \frac{\widehat{\gamma}_N}{K}, \quad \widehat{V}ar_{intra}\widehat{\mu}'_K = \frac{1}{K^2}\sum_{k\leq K}\frac{\widehat{V}_k}{n_k}.$$

Under a convergence rate for the sub-sample sizes, the estimators are all consistent, and the variances of order $K^{-1}$ are larger than the variance of $\widehat{\mu}_N$. The studentized statistics, as well as the normalized statistics converge in distribution:

**Proposition 2.1** *If $E|X|^{2+\delta} < \infty$ for some $\delta > 0$ and $n_k$ is of order $K^\alpha$, $0 < \alpha < 1/2$, for $k = 1, \ldots, K$, then $\widehat{\mu}_N$ converges a.s. to $\mu$ as $N \to \infty$, and $\widehat{\mu}'_K$ converges in probability to $\mu$ as $K \to \infty$, $\widehat{V}ar\widehat{\mu}_N - Var\widehat{\mu}_N$ and $\widehat{V}ar\widehat{\mu}'_K - Var\widehat{\mu}'_K$ converge in probability to 0, $(\widehat{V}ar\widehat{\mu}_N)^{-1/2}(\widehat{\mu}_N - \mu)$ and $(\widehat{V}ar\widehat{\mu}'_K)^{-1/2}(\widehat{\mu}'_K - \mu')$ converge in distribution to standard Gaussian variables.*

*Proof.* $\widehat{\mu}_N - \mu$ and $\widehat{\mu}'_K - \mu'$ are the means of the (weighted) independent variables having zero mean and variances of main order $n^*$ and $K$ respectively. With $n_k = O(K^\alpha)$ for $k = 1, \ldots, K$, $K = O(N^{1/(1+\alpha)})$, $n^* = O(N^{-1}K^\alpha) = O(K^{-1})$, and $K\widetilde{n} = O(K^{\alpha+1})$, the variances $S_N^2 = N^{-1}\sigma^2 + \gamma n^*$ of $\widehat{\mu}_N$ and $S_K^{'2} = K^{-1}\{\gamma + \sigma^2/\widetilde{n}\}$ of $\widehat{\mu}'_K$ ((7) and (8)) are such that $\gamma N^{-1}K^\alpha < S_N^2 \leq N^{-1}K^\alpha(\gamma + o(1))$ and $K^{-1}\gamma \leq S_K^{'2} \leq K^{-1}(\gamma + o(1))$. Then $\sum_N Var\widehat{\mu}_N < \infty$ which ensures the a.s. converge of the estimator $\widehat{\mu}_N$ by the Borel-Cantelli lemma.

The boundedness of the variances $S_N^2$ and $S_K^{'2}$ imply the following Lindeberg conditions

$$\frac{1}{S_N^2}\sum_{k,i} E(X_{ki}-\mu)^2 1_{\{|X_{ki}-\mu|>\varepsilon S_N\}} \stackrel{N\to\infty}{\longrightarrow} 0$$

$$\frac{1}{S_K^{'2}}\sum_k E(\widehat{\mu}_k-E\mu_k)^2 1_{\{|\widehat{\mu}_k-\mu|>\varepsilon S'_K\}} \stackrel{K\to\infty}{\longrightarrow} 0$$

and the weak convergence of the normalized variables is a consequence of a CLT [5]. ∎

If we consider $\bar{\mu}_K = K^{-1}\sum_{k=1}^K \mu_k$, similar results hold for $\widehat{\mu}'_K - \bar{\mu}_K$, with variance $\sigma^2(K\widetilde{n})^{-1} = Var_{intra}\widehat{\mu}'_K$, and for $\bar{\mu} - \mu$ with variance $\gamma K^{-1} = Var_{inter}\widehat{\mu}'_K$. The result in the proposition below is similar to the convergence property in a stratified population with fixed strata [4].

**Proposition 2.2** *Under the conditions of the proposition 2.1, $\widehat{\mu}'_K - \bar{\mu}_K$ and $\bar{\mu}_K - \mu'$ converge in probability to zero, $K^{1/2}(K^{-1}\sum_k n_k^{-1}\widehat{V}_k)^{-1/2}(\widehat{\mu}'_K - \bar{\mu})$ and $K^{1/2}\widehat{\gamma}_K^{-1/2}(\bar{\mu}_K - \mu')$ converge in distribution to standard Gaussian variables.*

## 3. Bootstrap estimation

The observed data set is denoted by $\mathcal{X}$ and we consider the three bootstrap procedures described in the introduction. For each of them, $E_*$ and $Var_*$ denote the mean and variance for the bootstrap sampling distribution conditionally on $\mathcal{X}$, without any distinction of the specific distribution, $\widehat{F}_k$ is the empirical sub-distribution function of the $k$-th population.

### 3.1. Bootstrap sampling of individuals

We consider a bootstrap sample which consists of the sub-samples set of size $n_k$ for the $k$-th observed population, $(X^*_{ki})_{i \leq n_k}$, where the variables $X^*_{ki}$ have the distribution $\widehat{F}_k$, $i \leq n_k, k \leq K$ and the $K$ populations are independent and considered as fixed. The bootstrap version of the previous estimators are written

$$\widehat{\mu}^*_k = = n_k^{-1} \sum_{i \leq n_k} X^*_{ki},$$

$$\widehat{\mu}^*_N = N^{-1} \sum_{k \leq K} \sum_{i \leq n_k} X^*_{ki} = N^{-1} \sum_{k \leq K} n_k \widehat{\mu}^*_k,$$

$$\widehat{\mu}'^*_K = K^{-1} \sum_{k \leq K} \widehat{\mu}^*_k,$$

they satisfy $E_* \widehat{\mu}^*_k = \widehat{\mu}_k$, $E_* \widehat{\mu}^*_N = \widehat{\mu}_N$, $E_* \widehat{\mu}'^*_K = \widehat{\mu}'_K$, $E_* \widehat{\mu}^{*2}_k = (n_k-1) n_k^{-2} \widehat{V}_k + \widehat{\mu}_k^2$. The variances of $\widehat{\mu}^*_N$ and $\widehat{\mu}'^*_K$ reduce to intra-population variances given by

$$Var_* \widehat{\mu}^*_k = \frac{n_k - 1}{n_k^2} \widehat{V}_k,$$

$$Var_* \widehat{\mu}^*_N = \frac{1}{N^2} \sum_{k \leq K} (n_k - 1) \widehat{V}_k,$$

$$Var_* \widehat{\mu}'*_K = \frac{1}{K^2} \sum_{k \leq K} \frac{n_k - 1}{n_k^2} \widehat{V}_k.$$

Finally, let

$$\widehat{V}^*_k = \frac{n_k}{(n_k - 1)^2} \sum_{i \leq n_k} (X^*_{ki} - \widehat{\mu}^*_k)^2,$$

it is an unbiased bootstrap estimator of $\widehat{V}_k$ and unbiased bootstrap estimators the variances of $\widehat{\mu}^*_N$ and $\widehat{\mu}'^*_K$ follow.

Then bootstrapping on independent individuals only with fixed populations (resampling scheme B2) as in the classical stratified sampling leads to estimators of $\mu$ having variances asymptotically equivalent to the within-populations variance of the estimators $\widehat{\mu}$. All variables $X^*_{ki}, i \leq n_k, k \leq n$ are independent conditionally on the sample $\mathcal{X}$ and $\widehat{\mu}^*_N - \widehat{\mu}_N$ is the mean of the $N$ independent centered variables $(X^*_{ki} - \widehat{\mu}_k)$. The next convergence results are then proved as Proposition 2.1,

**Proposition 3.1** *If $E|X|^{2+\delta} < \infty$ for some $\delta > 0$ and $n_k$ is of order $K^\alpha$, $0 < \alpha < 1/2$, for $k = 1, \ldots, K$, then conditionally on $\mathcal{X}$, $(\widehat{\mu}_N^* - \widehat{\mu}_N)$ and $(\widehat{\mu}_K'^* - \widehat{\mu}_K')$ converge a.s. to zero, $(\widehat{V}ar_{intra}^*\widehat{\mu}_N^*)^{-1/2}(\widehat{\mu}_N^* - \widehat{\mu}_N)$ and $(\widehat{V}ar_{intra}^*\widehat{\mu}_K'^*)^{-1/2}(\widehat{\mu}_K'^* - \widehat{\mu}_K')$ converge weakly to standard Gaussian variables.*

### 3.2. Bootstrap sampling of the populations

A sample of $K$ populations is drawn uniformly from the observed population set, each of them having the probability $K^{-1}$, and the variables $X_{ki}^*$ are defined as the observed values of the variable $X$ in these sampled populations, according to the sampling scheme B1. We get $\widehat{\mu}_k^* = \widehat{\mu}_l$ with probability $K^{-1}$, for every $k, l \leq K$, and $E_*\widehat{\mu}_k^* = \widehat{\mu}_K'$ and $E_*\widehat{\mu}_k^{*2} = K^{-1}\sum_{l \leq K}\widehat{\mu}_l^2$ for every $k$. It follows that both $\widehat{\mu}_K'^*$ and $\widehat{\mu}_N^*$ have the bootstrap mean $\widehat{\mu}_K'$ and the bootstrap variances are

$$Var_*\widehat{\mu}_N^* = N^{-1}\sum_{k \leq K} n_k(\widehat{\mu}_k - \widehat{\mu}_K')^2 = n^*\frac{K-1}{K}\widehat{V}ar\widehat{\mu}_K',$$

$$Var_*\widehat{\mu}_K'^* = \frac{1}{K^2}\sum_{k \leq K}(\widehat{\mu}_k - \widehat{\mu}_K')^2 = \frac{K-1}{K^2}\widehat{V}ar\widehat{\mu}_K'.$$

The usual bootstrap estimator of the variance of $\widehat{\mu}_K'^*$ is $K^{-1}\widehat{S}_K^{*2}$ where $\widehat{S}_K^{*2} = \frac{1}{K-1}\sum_{k \leq K}(\widehat{\mu}_k^* - \widehat{\mu}_K'^*)^2$ is the bootstrap variance of the $\widehat{\mu}_k^*$'s, it is a strongly consistent bootstrap estimator and $E^*\widehat{S}_K^{*2} = (K-1)K^{-1}\widehat{V}ar\widehat{\mu}_K'$. We now get bootstrap estimators of $\mu$ having a variance asymptotically equivalent to $K^{-1}\widehat{V}ar\widehat{\mu}_K'$ as $K \to \infty$ and the asymptotic behavior of $\widehat{\mu}_K'^*$ is similar to that of $\widehat{\mu}_K'$:

**Proposition 3.2** *Under the conditions of Proposition 3.4 and conditionally on $\mathcal{X}$, $K^{1/2}\widehat{S}_K^{*-1}(\widehat{\mu}_K'^* - \widehat{\mu}_K')$ converges weakly to a standard Gaussian variables.*

As the uniform population sampling is not relevant for $\widehat{\mu}_N$, let us consider a sample of $K$ populations drawn with probabilities $n_k N^{-1}$ for population $k = 1, \ldots, K$, and the variables $X_{ki}^*$ are defined as the observed values of the variable $X$ in these sampled populations. The means of the bootstrap estimators are $E_*\widehat{\mu}_k^* = E_*\widehat{\mu}_N^* = \widehat{\mu}_N$, and for $l = 1, \ldots, K$ their variance are

$$Var_*\widehat{\mu}_k^* = \frac{1}{N}\sum_{l \leq K} n_l(\widehat{\mu}_l - \widehat{\mu}_N)^2,$$

$$Var_*\widehat{\mu}_N^* = \frac{n^*}{N}\sum_{k \leq K} n_k(\widehat{\mu}_k - \widehat{\mu}_N)^2.$$

An unbiased bootstrap estimator of the variance of $\widehat{\mu}_N$ is also an estimator of this variance:

$$\widehat{S}_N^{*2} = \frac{n^*}{n^* - 1}\sum_{k \leq K} n_k(\widehat{\mu}_k^* - \widehat{\mu}_N^*)^2,$$

a weighted bootstrap variance of the $\widehat{\mu}_k^*$'s.

**Proposition 3.3** *Under the conditions of Proposition 3.4 and conditionally on $\mathcal{X}$, $\widehat{S}_N^{*-1}(\widehat{\mu}_N^* - \widehat{\mu}_N)$ converges weakly to a standard Gaussian variable.*

### 3.3. Cluster bootstrap sampling

For the B3 bootstrap sampling, $K$ independent populations are drawn uniformly from the observed population set and the bootstrap sample consists of sub-samples of size $n_l - 1$ and distribution $\widehat{F}_l$ if the the $k$-th bootstrap population is the $l$-th observed population. Denote $Y_{li}^*$ bootstrap variables having the distribution $\widehat{F}_l$ and $X_{ki}^*$ bootstrap variables having the distribution $\widehat{F}_l$ with probability $K^{-1}$, for every $l, k = 1, \ldots, K$. The bootstrap sampling distribution is then $\widehat{\widehat{F}} = K^{-1} \sum_{l \leq K} \widehat{F}_l$, $E_* X_{ki}^* = \widehat{\mu}_K'$ and $\widehat{\mu}_k^* = (n_l - 1)^{-1} \sum_{i \leq n_l} Y_{li}^*$ with probability $K^{-1}$, for every $l$ and $k = 1, \ldots, K$. We obtain $E_* \widehat{\mu}_k^* = E_* \widehat{\mu}_K'^* = \widehat{\mu}_K' = E_* \widehat{\mu}_N^*$,

$$E_* \widehat{\mu}_k^{*2} = \frac{1}{K} \sum_l \left\{ \frac{1}{n_l(n_l - 1)} \sum_{j \leq n_l} (X_{lj}^2 - \widehat{\mu}_l^2) + \widehat{\mu}_l^2 \right\} = \frac{1}{K} \sum_l \left\{ \frac{\widehat{V}_l}{n_l} + \widehat{\mu}_l^2 \right\}$$

for every $k$ and

$$Var_* \widehat{\mu}_K'^* = \frac{1}{K^2} \sum_k \{(\widehat{\mu}_k - \widehat{\mu}_K')^2 + \frac{\widehat{V}_k}{n_k}\} = \frac{K-1}{K^2} \widehat{Var} \widehat{\mu}_K' + \widehat{Var}_{intra} \widehat{\mu}_K'$$

is unbiasedly estimated by

$$\widehat{Var}^* \widehat{\mu}_K'^* = \frac{1}{K(K-1)} \sum_k (\widehat{\mu}_k^* - \widehat{\mu}_K'^*)^2.$$

The behavior of $\widehat{\mu}_N^*$ also depends on $\widehat{\mu}_K'^2$ in this procedure, with

$$Var_* \widehat{\mu}_N^* = (n^* - 1) K^2 Var_* \widehat{\mu}_K'^* + (\widehat{\mu}_K'^2 - \widehat{\mu}_N^2)$$

and $\widehat{\mu}_N^2$ cannot be simply estimated in this way. For $\widehat{\mu}_K'$, the variance of the bootstrap estimator $\widehat{\mu}_K'^*$ in the two-stage cluster procedure is equivalent to the sum of the total and within-populations variances of $\widehat{\mu}_N$ as $N$ and $K \to \infty$, it is therefore unsuitable, and we get

**Proposition 3.4** *Under the conditions of Proposition and conditionally on $\mathcal{X}$, $(Var_* \widehat{\mu}_K'^*)^{-1/2} (\widehat{\mu}_K'^* - \widehat{\mu}_K')$ converges weakly to a standard Gaussian variable.*

### 3.4. Second order asymptotics

It appears that among the studied bootstrap procedures, only the second resampling scheme of B2 provides consistent estimators for $\mu$ and the variances of the estimators. Further expansions prove that the bootstrap estimator of $\mu$ satisfy the classical properties of bootstrap estimators in the independent case [1, 2, 8, 10, 17]. The variances are denoted $S_N^2 = Var \widehat{\mu}_N$, $\widehat{S}_N^2 = \widehat{Var} \widehat{\mu}_N$ and $\widehat{S}_{*N}^2 = \widehat{Var}_* \widehat{\mu}_N^*$.

**Proposition 3.5** *If $E|X|^3 < \infty$ and $n_k = c_k K^\alpha$, for some positive and bounded $c_k$ and $0 < \alpha < 1/2$, for $k = 1, \ldots, K$, then the bootstrap estimators based on the population sampling procedure, with probability $n_k N^{-1}$ for population $k$, satisfy a.s.*

$$\limsup_{K \to \infty} \sup_x K^{1/2+2\alpha} |P\{S_N^{-1}(\widehat{\mu}_N - \mu) \leq x\}$$
$$- P_*\{\widehat{S}_N^{-1}(\widehat{\mu}_N^* - \widehat{\mu}_N) \leq x\}| \leq \frac{4CE|X-\mu|^3}{\gamma^{3/2}}.$$

*where $C$ is the constant of the Berry-Esseen bound.*

*If $E|X|^{6+\delta} < \infty$ for some $\delta > 0$ and the d.f. $F_k$ are continuous for $k \leq K$, then a.s.*

$$\limsup_{K \to \infty} \sup_x K^{1/2+2\alpha} |P\{\widehat{S}_N^{-1}(\widehat{\mu}_N - \mu) \leq x\} - P_*\{S_N^{*-1}(\widehat{\mu}_N^* - \widehat{\mu}_N) \leq x\}| = 0.$$

*Proof.* The variances satisfy $S_N^2 = K^{-1}(\gamma + \sigma^2 O(K^{-\alpha}))$, $\widehat{S}_N^2 = S_N^2 + o_p(1)$ a.s. under convergence rate of the $n_k$'s. The bootstrap variance $\widehat{S}_{*N}^2$ is estimated by $S_N^{*2}$ in 3.1, and $S_N^{*2} = S_N^2 + o(1)$ a.s. The Berry-Esseen theorem for a weighted sum of independent variables with varying distributions applies for both $\widehat{\mu}_N - \mu$ and $\widehat{\mu}_N^* - \widehat{\mu}_N$. Using the expansion $S_N^{-3} = K^{3/2}\{\gamma^{-3/2} + O(K^{-\alpha})\}$ and the inequality $E|N^{-1}n_k(\widehat{\mu}_k - \mu)|^3 \leq N^{-3}E|\sum_{i \leq n_k}(X_{ki} - \mu)|^3 \leq 4N^{-3}n_k E|X-\mu|^3$ for each $k$, we have $\sum_k E|N^{-1}n_k(\widehat{\mu}_k - \mu)|^3 \leq 4N^{-2}E|X-\mu|^3$ and, uniformly in $x$,

$$\sup_x |P\{S_N^{-1} \sum_{k \leq K} \frac{n_k}{N}(\widehat{\mu}_k - \mu) \leq x\} - \Phi(x)| \leq C S_N^{-3} \sum_k E(\frac{n_k}{N}|\widehat{\mu}_k - \mu|)^3$$
$$\leq K^{-1/2-2\alpha} \gamma^{-3/2} 4CE|X-\mu|^3 (1 + o_p(1)).$$

For the bootstrap estimator, a.s. conditionally on $\mathcal{X}$ uniformly in $x$,

$$E_* N^{-3} \sum_k n_k^3 |\widehat{\mu}_k^* - \widehat{\mu}_N|^3 = N^{-3} \sum_k n_k^3 |\widehat{\mu}_k - \mu|^3,$$

$$\sup_x |P_*\{\widehat{S}_N^{-1} \sum_{k \leq K} \frac{n_k}{N}(\widehat{\mu}_k^* - \widehat{\mu}_N) \leq x\} - \Phi(x)| \leq C\widehat{S}_N^{-3} \sum_k E_*|\frac{n_k}{N}(\widehat{\mu}_k^* - \widehat{\mu}_N)|^3$$
$$= CS_N^{-3} \sum_k (\frac{n_k}{N}|\widehat{\mu}_k - \mu|)^3 (1 + o_p(1))$$

and the first result follows.

The second part of the proposition is a consequence of Edgeworth expansions of $\widehat{\mu}_N - \mu$ and $\widehat{\mu}_N^* - \widehat{\mu}_N$, with $S_N^{-3} = K^{3/2}\{\gamma^{-3/2} + O(K^{-\alpha})\}$. Uniformly in $x$,

$$P\{\widehat{S}_N^{-1}(\widehat{\mu}_N - \mu) \leq x\} =$$
$$\Phi(x) + \frac{1}{6S_N^3} E(\sum_k |\frac{n_k}{N}(\widehat{\mu}_k - \mu)|^3) p(x)\phi(x) + o(K^{-(\alpha+1/2)}).$$

with $p(x) = 2x^2 + 1$. The expansion for the bootstrap estimator is

$$P_*\{(\widehat{S}_N^*)^{-1}(\widehat{\mu}_N^* - \widehat{\mu}_N) \leq x\} = \Phi(x) + \frac{1}{6\widehat{S}_N^{3/2}} E_*(\widehat{\mu}_k^* - \widehat{\mu}_N)^{3/2}\{p(x)\phi(x)$$

$$+ o(1)\} = \Phi(x) + \frac{1}{6\widehat{S}_N^{3/2}} N^{-3} \sum_k |n_k(\widehat{\mu}_k - \mu)|^3 \{p(x)\phi(x) + o(1)\},$$

with $\widehat{S}_N$ converging to $S_N$ and the sum of third moments converging to their expectation at the rate $N^{-1}$. ∎

## 4. Discussion

The model studied in this paper is the model of means of a variable $X_i$ having a mixture distribution of $L$ d.f. $F_1, \ldots, F_L$, defined as $F_k(x) = P\{X_i \leq x | i \in \mathcal{P}_k\}$, for the k-th sub-population $\mathcal{P}_k$. The variable $X_i$ has the d.f. $F = \sum_{k=1}^L p_k F_k$, with $p_k = P\{X_i \in \mathcal{P}_k\}$.

For a stratified population with a fixed number of strata, [2] have proved convergence of the bootstrap Studentized estimators under the condition of sub-sample sizes $n_k$ of the same order $O(N)$, they have proved in particular that the distribution function of $(Var\widehat{\mu}_k^*)^{-1/2} \sum_k E_*(\widehat{\mu}_k^* - \widehat{\mu}_k)^3$ is an approximation of $(Var\widehat{\mu}_k)^{-1/2} \sum_k E(\widehat{\mu}_k - \mu_k)^3$. Their model is different from the random model considered in this paper under the slower convergence rate $n_k = O(K^\alpha)$, $0 < \alpha < 1/2$ and with weighted estimators. In proposition 3.5, the order of the approximations is $K^{1/2+2\alpha}$, that is stronger than the usual results.

Inversion of several expansions improve coverage probabilities and have been compared to bootstrap confidence intervals deduced from the approximation of the distribution of Studentized statistics, these numerical improvements were discussed in [1, 8, 10] and in many other papers later. Here all the estimators are differentiable transformations of moment estimators that admit Edgeworth expansions which provides a second order approximation of the distribution function of $S_N^{-1}(\widehat{\mu}_N - \mu)$ by its inversion,

$$P\{S_N^{-1}(\widehat{\mu}_N - \mu) \leq x - (1 - x^2) \frac{\sum_k (\frac{n_k}{N} |\widehat{\mu}_k - \mu|)^3}{6\widehat{S}_N^3}\} = \Phi(x) + 0(K^{-(\alpha+1/2)})$$

but the order $0(K^{-(\alpha+1/2)})$, coming from the expression of $S_N$, differs from $O(K^{-1})$, the expected term after a correction term $O(K^{-1/2})$ in a one-stage sampling. For the studentized statistic, the correction is similar, with $1 - x^2$ replaced by $1 + 2x^2$.